\documentclass[11pt]{amsart}

\usepackage{amscd}
\usepackage{amsmath, amssymb}
\usepackage{amsfonts}

\newcommand{\de}{\partial}

\newcommand{\ov}[1]{\overline{#1}}

\newcommand{\ti}[1]{\tilde{#1}}

\newcommand{\vol}{\mathrm{Vol}}
\newcommand{\diam}{\mathrm{diam}}

\renewcommand{\leq}{\leqslant}
\renewcommand{\geq}{\geqslant}

\begin{document}
\newtheorem{theorem}{Theorem}[section]

\newtheorem{corollary}[theorem]{Corollary}
\newtheorem{proposition}[theorem]{Proposition}
\theoremstyle{definition}
\newtheorem{question}[theorem]{Question}
\theoremstyle{definition}
\newtheorem{example}[theorem]{Example}

\title{Calabi-Yau manifolds and their degenerations}
\author{Valentino Tosatti}
 \address{Department of Mathematics \\ Columbia University \\ New York, NY 10027}

  \email{tosatti@math.columbia.edu}
\begin{abstract}
Calabi-Yau manifolds are geometric objects of central importance in several branches of mathematics, including differential geometry, algebraic geometry and
mathematical physics. In this article we give a brief introduction to the subject
aimed at a general mathematical audience, and we present some of our results that shed some light on the possible ways in which families of Calabi-Yau manifolds can degenerate.
\end{abstract}
\thanks{This work was partially supported by National Science Foundation grant DMS-1005457.}

\maketitle
While preparing this article we faced a difficult choice: should the article be
addressed to a general science audience, or should it be written primarily for scholars in our field? We decided to strike a balance between these two approaches and write an article that can be read by any mathematician. We apologize to all the other potential readers for this choice. The interested reader can find a gentle introduction to this topic in the recent book of Yau and Nadis \cite{YN}
as well as in their article \cite{YN2}.

\section{Definitions of Calabi-Yau manifolds}

The main object of study in this article are {\em Calabi-Yau manifolds}. There are many possible definitions of these spaces, and we will start by reviewing a few of them. First of all recall that a complex manifold $X$ is a smooth manifold of real dimension $2n$ with an atlas whose transition functions are holomorphic. In particular each tangent space of $X$ is naturally identified with $\mathbb{C}^n$
and multiplication by $i$ induces a tensor $J:TX\to TX$ with $J^2=-{\rm Id}$. A Riemannian metric $g$ on $X$ is called Hermitian if $g(JY,JZ)=g(Y,Z)$ for all vectors $Y,Z$. In this case we define $\omega(Y,Z)=g(JY,Z)$, which is a skew-symmetric $2$-form on $X$. If $d\omega=0$ we say that $g$ is a K\"ahler metric, and we say that $X$ is K\"ahler if it admits K\"ahler metrics. The cohomology class $[\omega]$ lives in $H^2(X,\mathbb{R})\cap H^{1,1}_{\ov{\de}}(X)=:H^{1,1}(X,\mathbb{R}),$ and is called a K\"ahler class.
As we just explained, the tangent bundle of a complex manifold inherits a complex structure $J$, and so it has well-defined Chern classes $c_i(X)\in H^{2i}(X,\mathbb{Z})$, $1\leq i\leq n$.
We can now give several equivalent definitions of Calabi-Yau manifolds:

\begin{enumerate}
\item (Complex Geometry) A Calabi-Yau manifold is a compact K\"ahler manifold $X$ with first Chern class
$c_1(X)$ equal to zero in the cohomology group $H^2(X,\mathbb{R})$.

\item (Algebraic Geometry) A Calabi-Yau manifold is a compact K\"ahler manifold $X$ with torsion canonical bundle $K_X=\Lambda^n T^{1,0}X^*$, so that $\ell K_X\cong \mathcal{O}_X$ for some integer $\ell\geq 1$.

\item (Einstein Equation) A Calabi-Yau manifold is a compact complex manifold $X$ with a K\"ahler metric $\omega$ with Ricci curvature identically zero (Ricci-flat).

\item (Riemannian Geometry) A Calabi-Yau manifold is a compact Riemannian manifold $(X,g)$ of real dimension $2n$ with restricted holonomy group contained in the special unitary group $SU(n)$.

\end{enumerate}

We will explain why these definitions are equivalent after giving a few
examples.

\section{Examples of Calabi-Yau manifolds}
The following are some
simple examples of Calabi-Yau manifolds.\\

\begin{example} Let $X=\mathbb{C}^n/\Lambda$ be the quotient of Euclidean space $\mathbb{C}^n$ by a lattice
$\Lambda\cong\mathbb{Z}^{2n}$. Then $X$ is topologically a torus $(S^1)^{2n}$ and it has trivial tangent bundle and therefore also trivial canonical bundle. All Calabi-Yau manifolds of complex dimension $n=1$ are tori, and
are also called {\em elliptic curves}.
\end{example}

\begin{example} A Calabi-Yau manifold with complex dimension $n=2$ which is also simply connected is called a
{\em K3 surface}. Every Calabi-Yau surface is either a torus, a $K3$ surface, or a finite unramified
quotient of these. In general these quotients will have torsion but nontrivial canonical bundle, as is the case for
example for {\em Enriques surfaces} which are $\mathbb{Z}/2$ quotients of $K3$.
\end{example}

\begin{example} Let $X$ be a smooth complex hypersurface of degree $n+2$ inside complex projective space $\mathbb{CP}^{n+1}$.
Then by the adjuction formula the canonical bundle of $X$ is trivial, and so $X$ is a Calabi-Yau manifold. When $n=1$
we get an elliptic curve and when $n=2$ a $K3$ surface. More generally one can consider smooth
complete intersections in product of projective spaces, with suitable degrees, and get more examples of Calabi-Yau manifolds. \end{example}

\begin{example} Let $T=\mathbb{C}^2/\Lambda$ be a torus of complex dimension $2$ and consider the reflection through
the origin $i:\mathbb{C}^2\to \mathbb{C}^2$. This descends to an involution of $T$ with $16$ fixed points,
and we can take the quotient $Y=T/i$ which is an algebraic variety with $16$ singular rational double points (also known as
{\em orbifold} points). We resolve these $16$ points by blowing them all up and we get a map $f:X\to Y$ where $X$ is a smooth $K3$ surface,
known as the {\em Kummer surface} of the torus $T$.
\end{example}

\begin{example} A Calabi-Yau manifold $(X,g)$ of even complex dimension $n$ with holonomy equal to the symplectic group $Sp(\frac{n}{2})\subset SU(n)$ is called {\em hyperk\"ahler}.
We have that $Sp(1)=SU(2)$ so the only hyperk\"ahler manifolds of complex dimension $2$ are $K3$ surfaces. There are not many examples of higher-dimensional hyperk\"ahler manifolds. The simplest one is obtained by taking a $K3$ surface $Y$ and looking at $\ti{X}\to Y\times Y$ the blow-up of the diagonal. Flipping the two factors of $Y\times Y$ induces a $\mathbb{Z}/2$-action on $\ti{X}$ and the quotient space $X$ is hyperk\"ahler.
\end{example}

Let us now discuss why the four different definitions of Calabi-Yau manifolds that we gave are equivalent. The fact that definitions (3) and (4) are equivalent is
a simple exercise in Riemannian geometry, and to see that (2) implies (1) it suffices to take the first Chern class of the canonical bundle. The implication from (1) to (3) is the content of the celebrated Calabi Conjecture \cite{Ca}

\begin{theorem}[Yau's solution of the Calabi Conjecture \cite{yaupnas, yau1}]\label{yauthm}
On any compact K\"ahler manifold $X$ with $c_1(X)=0$ in $H^2(X,\mathbb{R})$ there exist K\"ahler metrics with Ricci curvature identically zero. Moreover, there is a unique such Ricci-flat metric in each
K\"ahler class of $X$.
\end{theorem}

Finally, the fact that (3) implies (2) is a consequence of a decomposition theorem due to Yau (cfr. \cite{Be, Ca}):
every compact Ricci-flat K\"ahler manifold of complex dimension $n$ has a finite unramified cover $\ti{X}$ which splits isometrically as a product of a flat torus, of simply connected Calabi-Yau manifolds with holonomy equal to $SU(n)$ and of simply connected hyperk\"ahler manifolds. In particular the canonical bundle of $\ti{X}$ is trivial,
which implies that the canonical bundle of $X$ is torsion.

Let us now briefly discuss why Calabi-Yau manifold are important in several branches of mathematics.
They are important in differential geometry because they give examples of Einstein metrics, which are Ricci-flat but not flat if $X$ is not covered by a torus, and these metrics are almost never explicit. Furthermore, if the K\"ahler class is fixed, the Ricci-flat metric is uniquely determined by the complex structure. Using this fact, Tian \cite{tian} and Todorov \cite{tod} proved that the moduli space of polarized Calabi-Yau manifolds is smooth.

They are important in algebraic geometry because of the position they occupy in the theory of classification
of algebraic varieties. These are usually divided into families according to their Kodaira dimension $\kappa(X)\in\{-\infty, 0,1,\dots,n\}$ where $n=\dim_{\mathbb{C}}(X)$, and conjecturally every algebraic variety with $\kappa(X)=0$ is birational to a Calabi-Yau variety (possibly singular). There are many other basic questions about Calabi-Yau manifolds which are still open: does every simply connected Calabi-Yau manifold have a rational curve?
Is the number of deformation classes of simply connected Calabi-Yau manifolds of a given dimension finite?

Calabi-Yau threefolds with holonomy $SU(3)$ are important in mathematical physics because they can be used
in string theory to construct supersymmetric theories which live (locally) on the product of four-dimensional spacetime with a Calabi-Yau threefold. This led to the discovery of the mathematical phenomenon of mirror symmetry,
which has generated a huge amount of research, see for example \cite{SYZ}. In the general framework of mirror symmetry one has a family of Calabi-Yau manifolds parametrized by a punctured disc $\Delta^*$ in the complex
plane, which degenerates when approaching the origin. The chief example of this is the family of Calabi-Yau quintic hypersurfaces $X_t$ in $\mathbb{CP}^4$ given by the equation in homogeneous coordinates
$$Z_0^5+Z_1^5+Z_2^5+Z_3^5+Z_4^5=5t^{-1}Z_0Z_1Z_2Z_3Z_4.$$
As $t$ approaches zero, $X_t$ degenerates to a union of $5$ hyperplanes. This is an example of a {\em large complex structure limit}. To such a family one then associates a ``mirror family'' of Calabi-Yau manifolds $\check{X}_t$ with fixed complex structure and varying K\"ahler class, which approaches a {\em large K\"ahler structure limit}. Then for $t$ close to zero, one can relate invariants of the complex structure of $X_t$ to invariants of the symplectic structure of $\check{X}_t$. For an introduction to this circle of ideas see \cite{GHJ}.

Finally let us note here that there are many possible generalizations of the concept of Calabi-Yau manifolds:
noncompact ones, singular ones, non-K\"ahler ones, symplectic ones, and so on. We will not delve here into
all these concepts and their precise definitions, but refer the reader to the survey article of Yau \cite{yausur}
instead.

\section{Degenerations of Calabi-Yau manifolds}
From the discussion in the previous section it is apparent that the study of degenerations of families of Calabi-Yau manifolds is an interesting topic with ramifications in several different branches of mathematics. In the rest of this article, we will study the behaviour of Ricci-flat K\"ahler metrics on a fixed Calabi-Yau manifold as their K\"ahler class degenerates. If the complex structure is allowed to vary, similar results were obtained by Ruan-Zhang \cite{rz} and Rong-Zhang \cite{rz2} building upon our work.

First of all let us identify the parameter space for K\"ahler classes on a compact K\"ahler manifold $X$. Recall that a K\"ahler class on $X$ is a cohomology class $\alpha$ in $H^{1,1}(X,\mathbb{R})$ which can be written as $\alpha=[\omega]$ for some K\"ahler metric $\omega$. The set of all K\"ahler classes is called the {\em K\"ahler cone} of $X$ and is
an open convex cone $\mathcal{K}_X\subset H^{1,1}(X,\mathbb{R}),$
which has the origin as its vertex.

\begin{question}\label{q1}
What is the behaviour of these Ricci-flat K\"ahler metrics when the class $\alpha$
degenerates to the boundary of the K\"ahler cone?
\end{question}

This question was posed by many people, including Yau \cite{yau2, yau3}, Wilson \cite{wilson} and McMullen \cite{ctm}. To get a feeling for what the K\"ahler cone and its boundary represent geometrically, we start with the following observation.
If $V\subset X$ is a complex subvariety of complex dimension $k>0$, then
it is well known (from the work of Lelong) that $V$ defines a homology class
$[V]$ in $H_{2k}(X,\mathbb{Z})$. Moreover if $[\omega]$ is a K\"ahler class, the pairing $\langle [V],[\omega]^{\smile k}\rangle$ equals
$$\int_V \omega^k=k! \vol(V,\omega)>0,$$
the volume of $V$ with respect to the K\"ahler metric $\omega$ (Wirtinger's Theorem). It follows that if a class $\alpha$ is on the boundary of $\mathcal{K}_X$ and if $V$ is any complex subvariety then the pairing $\langle [V],\alpha^{\smile k}\rangle$ is nonnegative, and moreover a theorem of Demailly-P\u{a}un \cite{dp} shows
that there must be subvarieties $V$ with pairing zero. Therefore as we approach the class $\alpha$ from inside $\mathcal{K}_X$, these subvarieties have volume that
goes to zero, and the Ricci-flat metrics must degenerate (in some way) along these subvarieties.

We now make Question \ref{q1} more precise. On a compact Calabi-Yau manifold $X$
fix a nonzero class $\alpha_0$ on the boundary of $\mathcal{K}_X$ and let $\{\alpha_t\}_{0\leq t\leq 1}$ be a smooth path of classes in $H^{1,1}(X,\mathbb{R})$ originating at $\alpha_0$ and with $\alpha_t\in \mathcal{K}_X$ for $t>0$. Call $\omega_t$ the unique Ricci-flat K\"ahler metric on $X$ cohomologous to $\alpha_t$ for $t>0$, which is produced by Theorem \ref{yauthm}.

\begin{question}\label{q2}
What is the behaviour of the Ricci-flat metrics $\omega_t$ when $t$ goes to zero?
\end{question}

Of course we could also consider sequences of classes instead of a path, and all we are going to say works equally well in that case.
Notice that we are not allowing the class $\alpha_t$ to go to infinity in $H^{1,1}(X,\mathbb{R})$ as it approaches $\de\mathcal{K}_X$. Because of this, we can prove the following basic fact, independently discovered by Zhang \cite{rz}:

\begin{theorem}[Tosatti \cite{deg}, Zhang \cite{rz}] The diameter of the metrics $\omega_t$ has a uniform upper bound as $t$ approaches zero,
\begin{equation}\label{diameter}
\diam(X,\omega_t)\leq C.
\end{equation}
\end{theorem}

On the other hand it is easy to construct examples of Ricci-flat K\"ahler metrics with unbounded cohomology class that violate \eqref{diameter}, by just rescaling a fixed metric by a large number.

Going back to Question \ref{q2}, the problem splits naturally into two cases which exhibit a rather different behaviour, according to whether the total integral $\int_X \alpha_0^n$ is strictly positive or zero. If $\int_X \alpha_0^n$ is positive this means that the volume
$$n! \vol(X,\omega_t)=\int_X\omega_t^n=\int_X\alpha_t^n$$
remains bounded away from zero as $t\to 0$, and this is called the {\em non-collapsing} case. If $\int_X \alpha_0^n=0$ then the volume $\vol(X,\omega_t)$
converges to zero, and this is called the {\em collapsing} case.

The main Question \ref{q2} falls into the general problem of understanding limits of sequences of Einstein manifolds with an upper bound for the diameter (but no bound for the sectional curvature in general), a topic that has been extensively studied (see e.g. \cite{gr, and, bkn, CC}). Our results are of a quite different nature from these works
and give stronger conclusions. The first theorem gives a satisfactory answer in the noncollapsing case:

\begin{theorem}[Tosatti \cite{deg}]\label{degen} If $\int_X \alpha_0^n>0$ then the Ricci-flat metrics $\omega_t$ converge smoothly away from an analytic subvariety $S$ to an incomplete Ricci-flat metric on its complement.
\end{theorem}

In fact the subvariety $S$ in this theorem is simply the union of all complex subvarieties where $\alpha_0$ integrates to zero. Whenever $\alpha_0$ is a rational class, the limit can also be understood using algebraic geometry: the subvariety $S$ can be contracted to create a singular Calabi-Yau manifold, and the limit metric is the pullback of a Ricci-flat metric on the smooth part.\\

In the collapsing case, when $\int_X\alpha_0^n=0$, things are more complicated. If $X$ is projective and the class $\alpha_0$ is rational, then the Log Abundance Conjecture in algebraic geometry implies that there is a fibration
$f:X\to Y$ where $Y$ is an algebraic variety of strictly lower dimension $m$ and $\alpha_0$ is the
pullback of an ample class on $Y$. If we call $S$ the critical locus of $f$ inside $X$ then $S$ is a subvariety and $f:X\backslash S\to Y\backslash f(S)$ is a smooth submersion with fibers
Calabi-Yau manifolds $X_y=f^{-1}(y)$ of complex dimension $n-m$. The subvariety $S$ is the union of all singular fibers of $f$ together with all the fibers with dimensions strictly larger than $n-m$. We also take $\alpha_t=\alpha_0+t[\omega_X]$, where $\omega_X$ is a fixed K\"ahler metric on $X$.
We then have the following result,
which says that the Ricci-flat metrics shrink the manifold to the base of the fibration:

\begin{theorem}[Tosatti \cite{deg2}]\label{degen2} Let $f:X\to Y$ be such a holomorphic fibration
and $\alpha_t=\alpha_0+t[\omega_X]$.
Then there is a smooth K\"ahler metric $\omega$ on $Y\backslash f(S)$ such that
when $t$ approaches zero the Ricci-flat metrics $\omega_t$ converge to $f^*\omega$.
The metric $\omega$ has Ricci curvature equal to a Weil-Petersson metric that measures the change of complex structures of the Calabi-Yau fibers $X_y$.
\end{theorem}

If we furthemore assume that $X$ is projective and all the smooth fibers $X_y$ are tori, then we can prove a stronger result:
\begin{theorem}[Gross, Tosatti, Zhang \cite{gtz}]\label{degen3} In the same setting as Theorem \ref{degen2}, assume that
$X$ is projective and the smooth fibers $X_y$ are tori. Then the convergence of $\omega_t$ to $f^*\omega$
is smooth and the sectional curvature of $\omega_t$ remains locally
bounded on $X\backslash S$. Along each torus fiber $X_y$ the rescaled metrics $t^{-1}\omega_t|_{X_y}$ converge smoothly to a flat metric. Finally, for any Gromov-Hausdorff limit space $(\hat{X},d)$ of $(X,\omega_t)$ there
is a local isometric embedding $(Y\backslash S,\omega)\hookrightarrow (\hat{X},d)$ with dense image.
\end{theorem}

Theorem \ref{degen3} can also be applied to study the large complex structure limits of families of polarized hyperk\"ahler manifolds in the large complex structure limit, see \cite{gtz} for details.

\section{Examples of degenerations}\label{ex}
First of all notice that Question \ref{q2} is only interesting if
$\dim H^{1,1}(X,\mathbb{R})>1,$ because otherwise $\mathcal{K}_X$ reduces to
an open half-line and there is only one Ricci-flat K\"ahler metric on $X$ up
to global scaling by a constant, so the only possible degenerations are given by scaling this metric to zero or infinity. For this reason, Question \ref{q2} is essentially void on Calabi-Yau manifolds of dimension $n=1$ (i.e. elliptic curves).

\begin{example}\label{ex1} Let $X=\mathbb{C}^n/\Lambda$ be a complex torus. A Ricci-flat K\"ahler metric on $X$ is the same as a flat K\"ahler metric, and each flat metric can be identified simply with a positive definite Hermitian $n\times n$ matrix. The boundary of the K\"ahler cone is then represented by non-negative definite Hermitian matrices $H$ with nontrivial kernel $\Sigma\subset\mathbb{C}^n$ (notice that in this case every class on $\de \mathcal{K}_X$ has zero integral, so we are always in the collapsing case).

If the class $\alpha_0$ corresponds to such a matrix $H$ with the kernel $\Sigma$ which is $\mathbb{Q}$-defined modulo $\Lambda$, then we can quotient $\Sigma$ out and get a map $f:X\to Y=\mathbb{C}^m/\Lambda'$ to a lower-dimensional torus ($m<n$) such that $H=f^*H'$ with $H'$ a positive definite $m\times m$ Hermitian matrix. It follows that when $t$ approaches zero, the (Ricci-)flat metrics $\omega_t$ collapse to the flat metric on $Y$ that corresponds to $H'$.
This is a special case of Theorems \ref{degen2} and \ref{degen3}, with $S$ empty and Weil-Petersson metric identically zero (since all the torus fibers are isomorphic).

If on the other hand the kernel $\Sigma$ is not $\mathbb{Q}$-defined, then $\Sigma$ defines a foliation on $X$ (which is not a fibration anymore) and the limit $H$ of the (Ricci-)flat metrics is a smooth nonnegative form which is {\em transversal} to the foliation (that means, positive in the complementary directions).
This case gives an idea of what to expect in general in the collapsing case when there is no fibration.
\end{example}

\begin{example} Let $f:X\to Y$ be the Kummer $K3$ surface of a torus $T$,
where $Y=T/i$ is the singular quotient of $T$ and $f$ is the blowup map. Take $\alpha_0$ to be the pullback of an ample divisor on $Y$, and note that
$\int_X\alpha_0^2>0$. If we call $S$ the union of the $16$ exceptional divisors of $f$, that is the union of the $16$ spheres $S^2$ which are the preimages of the singular points of $Y$, then $S$ is a complex submanifold of $X$. Then Kobayashi-Todorov \cite{kt} (using classical results on the moduli space of $K3$ surfaces, such as the Torelli theorem) proved that for any path $\alpha_t$ of K\"ahler classes that approach $\alpha_0$, the Ricci-flat metrics $\omega_t$ converge smoothly away from $S$ to the pullback of
the unique flat orbifold metric on $Y$ cohomologous to the ample divisor we chose. Here an orbifold flat metric on $Y$ simply means a flat metric on $T$ which is invariant under $i$. This statement is a special case of our Theorem \ref{degen}, which in particular gives a new proof of the result of Kobayashi-Todorov.
\end{example}

\begin{example} Let $X$ be a $K3$ surface which admits an elliptic fibration $f:X\to\mathbb{CP}^1=Y$. This means that $f$ is a surjective holomorphic map with all the fibers smooth elliptic curves except a finite number of fibers $S$ which are singular elliptic curves. Again we take $\alpha_0$ to be the pullback of an ample divisor on $Y$ and note that $\int_X\alpha_0^2=0$. We also take $\alpha_t=\alpha_0+t[\omega_X]$ for a K\"ahler metric $\omega_X$.
 Then Gross-Wilson \cite{gw} have shown that when $t$ goes to zero the metrics $\omega_t$ converge smoothly away from $S$ to the pullback $f^*\omega$, where $\omega$ is a K\"ahler metric on $Y=\mathbb{CP}^1$ minus the finitely many points $f(S)$ with singular preimage. Moreover they showed that away from $S$
 the rescaled Ricci-flat metrics along the fibers $t^{-1}\omega_t|_{X_y}$ converge
to flat metrics. More recently Song-Tian \cite{st} have noticed that the metric $\omega$ on $\mathbb{CP}^1\backslash f(S)$ has Ricci curvature equal to the pullback of the Weil-Petersson metric from the moduli space of elliptic curves via the map that to a point in $\mathbb{CP}^1\backslash f(S)$ associates the elliptic curve which lies above that point. This result is a special case of our Theorems \ref{degen2} and \ref{degen3}, which also provide a new proof of the theorem of Gross-Wilson.
\end{example}

\section{Questions}

Let us now mention a few open problems related to the above results, which we find very interesting.

\begin{question}\label{qq1} In the same setting as Theorem \ref{degen2} prove that
the rescaled Ricci-flat metrics along the fibers $t^{-1}\omega_t|_{X_y}$ converge smoothly to
the unique Ricci-flat K\"ahler metric on $X_y$ cohomologous to $[\omega_X]|_{X_y}$.
\end{question}
As explained in \cite{gtz} this would be implied by the following:
\begin{question}\label{qq2} In the same setting as Theorem \ref{degen2} prove that
the convergence of $\omega_t$ to $f^*\omega$ is smooth away from $S$.
\end{question}

Both Question \ref{qq1} and \ref{qq2} are solved in \cite{gtz} when $X$ is projective and the smooth fibers are tori.

The natural remaining question is what happens to the Ricci-flat metrics $\omega_t$ when $\alpha_0$ is an irrational class with $\int_X\alpha_0^n=0$, so that
there is no fibration structure. We conjecture the following:

\begin{question}\label{qq3} In this situation there is a subvariety $S\subset X$ and a smooth nonnegative
$(1,1)$-form $\omega_0$ on $X\backslash S$, which satisfies $\omega_0^n=0$, so that the Ricci-flat metrics $\omega_t$ converge smoothly away from $S$ to $\omega_0$.
\end{question}

In this case taking the kernel of $\omega_0$ we would get a foliation on $X\backslash S$ with leaves holomorphic subvarieties, a so-called {\em Monge-Amp\`ere foliation}. Of course Question \ref{qq3} is correct in the case of tori, as in Example \ref{ex1}.

Another problem that seems very interesting is whether the convergence in Theorem \ref{degen2} holds in the Gromov-Hausdorff sense.
More precisely, consider the metric space completion $(Z,d)$ of $(Y\backslash f(S),\omega)$.

\begin{question}\label{qq4} In the setting of Theorem \ref{degen2}, do the Ricci-flat manifolds $(X,\omega_t)$ converge to $(Z,d)$ in the Gromov-Hausdorff sense? Is $Z$ homeomorphic to $Y,$ the algebro-geometric limit?
\end{question}

The Gromov-Hausdorff converge is proved for $K3$ surfaces in the noncollapsing case in \cite{deg}, and in the general setting of Theorem \ref{degen} it is a recent result of Rong-Zhang \cite{rz2}. The homeomorphism question remains open even in this case. In the case of collapsing $K3$ surfaces, question \ref{qq4} was answered positively by Gross-Wilson \cite{gw}.
Finally, we pose the following:
\begin{question} Let $X$ be a compact Calabi-Yau manifold and $\alpha$ a cohomology class on $\de\mathcal{K}_X$.
Does there exist a smooth nonnegative form $\omega$ cohomologous to $\alpha$?
\end{question}
If $X$ is projective and $\alpha$ is rational then this can be proved using algebraic geometry if
$\int_X\alpha^n>0$, and it follows from the Log Abundance Conjecture if $\int_X\alpha^n=0$, see \cite{deg}.


\begin{thebibliography}{99}
\bibitem{and} Anderson, M.T. \emph{Ricci curvature bounds and Einstein metrics on compact manifolds},
J. Amer. Math. Soc. {\bf 2} (1989), no. 3, 455--490.
\bibitem{bkn} Bando, S., Kasue, A., Nakajima, H. \emph{On a construction of coordinates at infinity on manifolds with fast curvature decay and maximal volume growth}, Invent. Math. {\bf 97} (1989), no. 2, 313--349.
\bibitem{Be} Beauville, A. {\em Vari\'et\'es K\"ahleriennes dont la premi\`ere classe de Chern est nulle}, J. Differential Geom. {\bf 18} (1983), no. 4, 755--782.
\bibitem{Ca} Calabi, E. {\em On K\"ahler manifolds with vanishing canonical class}, in {\em Algebraic geometry and topology. A symposium in honor of S. Lefschetz,} pp. 78--89. Princeton University Press, Princeton, N. J., 1957.
\bibitem{CC} Cheeger, J., Colding, T. {\em On the structure of spaces with Ricci curvature bounded below. I},  J. Differential Geom. {\bf 46}  (1997),  no. 3, 406--480.
\bibitem{dp} Demailly, J.-P., P\u{a}un, M. \emph{Numerical characterization of the K\"ahler cone of a compact K\"ahler manifold}, Ann. of Math. {\bf 159} (2004), no. 3, 1247--1274.
\bibitem{gr} Gromov, M. \emph{Metric structures for Riemannian and non-Riemannian spaces}, Birkh\"auser, Boston, 1999.
\bibitem{GHJ} Gross, M., Huybrechts, D., Joyce, D. \emph{Calabi-Yau
manifolds and related geometries,} Springer-Verlag 2003.
\bibitem{gtz} Gross, M., Tosatti, V., Zhang, Y. {\em Collapsing of Abelian fibred Calabi-Yau manifolds}, arXiv:1108.0967.
\bibitem{gw} Gross, M., Wilson, P.M.H. \emph{Large complex structure limits of $K3$ surfaces}, J. Differential Geom. {\bf 55} (2000), no. 3, 475--546.
\bibitem{kt} Kobayashi, R., Todorov, A.N. \emph{Polarized period map for generalized $K3$ surfaces and the moduli of Einstein metrics}, Tohoku Math. J. {\bf 39} (1987), no. 3, 341--363.
\bibitem{ctm} McMullen, C.T. \emph{Dynamics on $K3$ surfaces: Salem numbers and Siegel disks}, J. Reine Angew. Math. {\bf 545} (2002), 201--233.
\bibitem{rz2} Rong, X. Zhang, Y. {\em Continuity of extremal transitions and flops for Calabi-Yau manifolds}, arXiv:1012.2940.
\bibitem{rz} Ruan, W.-D., Zhang, Y. \emph{Convergence of Calabi-Yau manifolds}, Adv. Math. {\bf 228} (2011), 1543--1589.
\bibitem{st} Song, J., Tian, G. \emph{The K\"ahler-Ricci flow on surfaces of positive Kodaira dimension}, Invent. Math. {\bf 170} (2007), no. 3, 609--653.
\bibitem{SYZ} Strominger, A., Yau, S.-T., Zaslow, E. {\em Mirror symmetry is $T$-duality}, Nuclear Phys. B {\bf 479} (1996), no. 1-2, 243--259.
\bibitem{tian} Tian, G. \emph{Smoothness of the universal deformation space of compact Calabi-Yau manifolds and its Petersson-Weil metric}, in {\em Mathematical aspects of string theory (San Diego, Calif., 1986),} 629--646, Adv. Ser. Math. Phys., {\bf 1}, World Sci. Publishing, Singapore, 1987.
\bibitem{tod} Todorov, A.N. {\em The Weil-Petersson geometry of the moduli space of ${\rm SU}(n\geq 3)$ (Calabi-Yau) manifolds. I}, Comm. Math. Phys. {\bf 126} (1989), no. 2, 325--346.
\bibitem{deg} Tosatti, V. \emph{Limits of Calabi-Yau metrics when the K\"ahler class degenerates}, J. Eur. Math. Soc. (JEMS) {\bf 11} (2009), no.4, 755-776.
\bibitem{deg2} Tosatti, V. \emph{Adiabatic limits of Ricci-flat K\"ahler metrics},  J. Differential Geom. {\bf 84} (2010), no.2, 427--453.
\bibitem{wilson} Wilson, P.M.H. \emph{Metric limits of Calabi-Yau manifolds}, in \emph{The Fano Conference}, 793--804, Univ. Torino, Turin, 2004.
\bibitem{yaupnas} Yau, S.-T. \emph{Calabi's conjecture and some new results in algebraic geometry}, Proc. Nat. Acad. Sci. U.S.A. {\bf 74} (1977), no. 5, 1798--1799.
\bibitem{yau1} Yau, S.-T. \emph{On the Ricci curvature of a compact K\"ahler manifold and the complex Monge-Amp\`ere equation, I}, Comm. Pure Appl. Math. {\bf 31} (1978), 339--411.
\bibitem{yau2} Yau, S.-T. \emph{Problem section} in \emph{Seminar on Differential Geometry}, pp. 669--706,
Ann. of Math. Stud. {\bf 102}, Princeton Univ. Press, 1982 (problem 49).
\bibitem{yau3} Yau, S.-T. \emph{Open problems in geometry}, Proc. Sympos. Pure
Math. {\bf 54} (1993), 1--28 (problem 88).
\bibitem{yausur} Yau, S.-T. {\em A survey of Calabi-Yau manifolds}, in {\em Surveys in differential geometry. Vol. XIII. Geometry, analysis, and algebraic geometry: forty years of the Journal of Differential Geometry,} 277--318, Surv. Differ. Geom., {\bf 13}, Int. Press, Somerville, MA, 2009.
\bibitem{YN} Yau, S.-T., Nadis, S. {\em The shape of inner space.
String theory and the geometry of the universe's hidden dimensions}, Basic Books, New York, 2010.
\bibitem{YN2} Yau, S.-T., Nadis, S. {\em String theory and the geometry
of the Universe's hidden dimensions}, Notices Amer. Math. Soc. {\bf 58} (2011), no. 8, 1067--1076.
\end{thebibliography}
\end{document}